\documentclass[12pt,a4paper]{article}

\usepackage{amssymb}
\usepackage{amsmath}
\usepackage{amsthm}
\usepackage{xypic}
\theoremstyle{plain}
\newtheorem{theorem}{Theorem}
\newtheorem{proposition}[theorem]{Proposition}

\theoremstyle{remark}
\newtheorem{remark}[theorem]{Remark}
\theoremstyle{definition}

\newtheorem{lemma}[theorem]{Lemma}

\def\varinjlim_#1{\lim\limits_{\longrightarrow\atop{#1}}}

\def\id{\mathop{\rm id}\nolimits}

\begin{document}
\author{A. V. Ershov}
\title{Logarithms of formal groups \\ over Hopf algebras}
\date{}
\maketitle This paper is the continuation of section 2 of the
paper "Floating bundles and their applications" \cite{pap1}.

Let $(H,\mu,\eta,\Delta,\varepsilon,S)$
be a commutative Hopf algebra over ring $R$ without torsion
and ${\frak F}(x\otimes 1,1\otimes x)$ be a
formal group over Hopf algebra $H.$
By $H_\mathbb{Q}$ denote the Hopf algebra
$H{\mathop{\otimes}\limits_\mathbb{Z}}\mathbb{Q}$
over ring
$R_\mathbb{Q}=R{\mathop{\otimes}\limits_\mathbb{Z}}\mathbb{Q}.$
We shall write $\mu,\eta,\ldots$ instead of $\mu_\mathbb{Q},
\eta_\mathbb{Q},\ldots.$

The aim of this paper is to prove the following result.
\begin{proposition}
For any commutative formal group ${\frak F}(x\otimes 1,1\otimes x),$
which is considered as a formal group over $H_\mathbb{Q},$
there exists a homomorphism  to
a formal group of the form ${\frak c}+x\otimes 1+1\otimes x,$
where $\frak c\in H_\mathbb{Q}{\mathop{\widehat{\otimes}}
\limits_{R_\mathbb{Q}}}H_\mathbb{Q}$
such that $(\id \otimes \varepsilon){\frak c}=0=
(\varepsilon \otimes \id){\frak c}.$
\end{proposition}
We recall that the notion of a homomorphism of formal
groups over Hopf algebra was given in \cite{pap2}.
Below we shall use notations of paper \cite{pap1}.

To prove the Proposition, we need the following Lemma.
\begin{lemma}
A symmetric series of the form ${\frak c}+x\otimes 1+1\otimes x
\in H_\mathbb{Q}{\mathop{\widehat{\otimes}}
\limits_{R_\mathbb{Q}}}H_\mathbb{Q}[[x\otimes 1,1\otimes x]]$
is a formal group over the Hopf algebra $H_\mathbb{Q}$
if and only if the following two conditions hold:
\begin{itemize}
\item[(i)] $(\id \otimes \Delta){\frak c}+1\otimes
{\frak c}-(\Delta \otimes \id){\frak c}-{\frak c}\otimes 1=0;$
\item[(ii)] $(\id \otimes \varepsilon){\frak c}=0=
(\varepsilon \otimes \id){\frak c}.$
\end{itemize}
\end{lemma}
{\raggedright(Note} that the condition (i) means, that ${\frak c}$ is a
2-cocykle in the cobar complex of the Hopf algebra
$H_\mathbb{Q}.$)\\
{\raggedright {\it Proof of the Lemma}.}\;
The conditions (i) and (ii) are equivalent to the associativity
axiom and to the unit axiom for formal groups respectively.
Let us show that the series $\Theta(x)=-(\mu\circ (\id \otimes S))
{\frak c}-x$ is the inverse element. Indeed,
$$({\mu}\circ (\id \otimes S)){\frak c}+
x+\Theta(x)=
(\mu \circ (\id \otimes S)){\frak c}+x-
(\mu \circ (\id \otimes S)){\frak c}-x=0.$$
The symmetric condition follows from the equality
$({\mu}\circ (\id \otimes S)){\frak
c}=({\mu}\circ (S\otimes \id)){\frak c}.\quad \square$

{\raggedright {\it Proof of the Proposition}.}\;
By definition, put
$$\widetilde{\omega}(x)=(\id \otimes \widetilde{\varepsilon})
\frac{\partial{\frak F}(x,z)}{\partial z}\in H[[x]]$$
(here $\widetilde{\varepsilon}\colon H[[z]]\rightarrow R$
is the map such that $\widetilde{\varepsilon}\mid_H=\varepsilon
\colon H\rightarrow R,\; \widetilde{\varepsilon}(z)=0$).
Recall (\cite{pap1}) that $\widetilde{\Delta}\colon
H[[x]]\rightarrow H[[x]]{\mathop{\widehat{\otimes}}\limits_R}H[[x]]=
H{\mathop{\widehat{\otimes}}\limits_R}H[[x\otimes 1,1\otimes x]]$
is the map such that $\widetilde{\Delta}\mid_H=\Delta,\;
\widetilde{\Delta}(x)={\frak F}(x\otimes 1,1\otimes x).$
We have
$$
(\Delta \widetilde{\omega})({\frak F}(x\otimes 1,1\otimes x))=
\widetilde{\Delta}(\widetilde{\omega}(x))=(\widetilde{\Delta}
\circ (\id \otimes \widetilde{\varepsilon})\circ \frac{\partial}
{\partial z})({\frak F}(x,z))=
$$
$$
((\id \otimes \id \otimes \widetilde{\varepsilon})\circ
(\widetilde{\Delta}\otimes \id)\circ \frac{\partial}{\partial z})
({\frak F}(x,z))=
$$
$$
((\id \otimes \id \otimes \widetilde{\varepsilon})\circ
\frac{\partial}{\partial z}\circ (\widetilde{\Delta}\otimes
\id))({\frak F}(x,z))=
$$
$$
\left((\id \otimes \id \otimes \widetilde{\varepsilon})\circ
\frac{\partial}{\partial z}\right)((\Delta \otimes \id){\frak F})({\frak
F}(x\otimes 1,1\otimes x),z)=
$$
$$
\left((\id \otimes \id \otimes \widetilde{\varepsilon})\circ
\frac{\partial}{\partial z}\right)((\id \otimes \Delta){\frak F})(x
\otimes 1,{\frak F}(1\otimes x,z))=
$$
$$
(\id \otimes \id \otimes \widetilde{\varepsilon})\left(\frac{\partial
((\id \otimes \Delta){\frak F})(x\otimes 1,{\frak F}(1\otimes x,
z))}{\partial {\frak F}(1\otimes x,z)}\right)\cdot 1\otimes\left((\id
\otimes \widetilde{\varepsilon})\frac{\partial{\frak F}
(1\otimes x,z)}{\partial z}\right)=
$$
$$
\frac{\partial{\frak F}(x\otimes 1,1\otimes x)}
{\partial (1\otimes x)}\cdot(1\otimes \widetilde{\omega})(1\otimes x).
$$
Therefore, we have
\begin{equation}
\label{1}
(\Delta \widetilde{\omega})({\frak F}(x\otimes 1,1\otimes x))=
\frac{\partial{\frak F}(x\otimes 1,1\otimes x)}
{\partial (1\otimes x)}\cdot(1\otimes \widetilde{\omega})(1\otimes x).
\end{equation}
If $${\frak F}(x,z)=\sum_{i,j\geq 0}A_{i,j}x^iz^j\quad (A_{i,j}\in
H{\mathop{\widehat{\otimes}}\limits_R}H),$$
then
$$\widetilde{\omega}(x)=(\id \otimes \widetilde{\varepsilon})
\sum_{i,j}A_{i,j}x^ijz^{j-1}=(\id \otimes \varepsilon)A_{0,1}+
\sum_{i\geq 1}((\id \otimes \varepsilon)A_{i,1})x^i,
$$
where $(\varepsilon \circ(\id \otimes \varepsilon))A_{0,1}=1\neq 0.$
Therefore
$$\frac{1}{\widetilde{\omega}(x)}\in H[[x]]\quad \hbox{and}
$$
$$
\widetilde{\Delta}\left(\frac{1}{\widetilde{\omega}(x)}\right)=
\frac{1}{\widetilde{\Delta}(\widetilde{\omega}(x))}=
\frac{1}{(\Delta\widetilde{\omega})({\frak F}(x\otimes 1,1\otimes x))}
\in H{\mathop{\widehat{\otimes}}\limits_R}H[[x\otimes 1,1\otimes x]].
$$
Therefore (\ref{1}) may be rewritten in the form
\begin{equation}
\label{2}
\frac{d(1\otimes x)}{(1\otimes
\widetilde{\omega})(1\otimes x)}=\frac{d{\frak F}
(x\otimes 1,1\otimes x)}{(\Delta\widetilde{\omega})({\frak F}
(x\otimes 1,1\otimes x))}.
\end{equation}
It is clear that
$$
\frac{1}{\widetilde{\omega}(x)}=b_0+
b_1x+\ldots ,
$$
where $b_i\in H,\; \varepsilon(b_0)=1.$
By ${\frak g}(x)$ denote the series
$$\int_o^x\frac{dt}{\widetilde{\omega}(t)}
\in H_\mathbb{Q}[[x]].$$
Equality (\ref{2}) implies
\begin{equation}
\label{3}
{\frak c}'+(1\otimes {\frak g})(1\otimes x)=
(\Delta{\frak g})({\frak F}(x\otimes 1,1\otimes x)),
\end{equation}
where ${\frak c}'$ is independent of $1\otimes x.$
The application of $\id\otimes \widetilde{\varepsilon}$
to equation (\ref{3}) yields
$$
(\id \otimes \widetilde{\varepsilon}){\frak c}'=(((\id \otimes
\varepsilon)\circ \Delta){\frak g})(x\otimes 1)=({\frak g}\otimes
1)(x\otimes 1),
$$
and the application of $\; \widetilde{\varepsilon}\otimes \id \;$
to equation (\ref{3}) yields
$$
(\widetilde{\varepsilon}\otimes \id){\frak c}'
+(1\otimes{\frak g})(1\otimes x)=(1\otimes{\frak g})(1\otimes x).
$$
Hence
\begin{equation}
\label{4}
(\Delta{\frak g})({\frak F}(x\otimes 1,1\otimes x))={\frak c}+
({\frak g}\otimes 1)(x\otimes 1)+(1\otimes {\frak g})(1\otimes x),
\end{equation}
where ${\frak c}'=({\frak g}\otimes 1)(x\otimes 1)+{\frak c},\;
{\frak c}\in H_\mathbb{Q}{\mathop{\widehat{\otimes}}
\limits_{R_\mathbb{Q}}}H_\mathbb{Q}\; \; \hbox{and}\; \;
(\id \otimes \varepsilon)
{\frak c}=0=(\varepsilon \otimes \id){\frak c}.$

To complete the proof we must check the condition (i)
of the previous Lemma. For this purpose we apply
$\id \otimes \widetilde{\Delta}\; \hbox{and}\; \widetilde{\Delta}
\otimes \id$ to equation ({\ref{4}}).
We have
$$
(((\id \otimes \Delta)\circ \Delta){\frak g})
((\id \otimes \Delta){\frak F}(x
\otimes 1\otimes 1,1\otimes {\frak F}(x\otimes 1,1\otimes x)))=
$$
$$
(\id \otimes \Delta){\frak c}+{\frak g}(x)\otimes 1\otimes 1+
1\otimes (\Delta {\frak g})({\frak F}(x\otimes 1,1\otimes x))=
$$
$$
(((\Delta \otimes \id)\circ \Delta){\frak g})
((\Delta \otimes \id){\frak F}({\frak F}(x\otimes 1,1\otimes x)
\otimes 1,1\otimes 1 \otimes x))=
$$
$$
(\Delta \otimes \id){\frak c}+(\Delta {\frak g})
({\frak F}(x\otimes 1,1\otimes x))\otimes 1+1\otimes 1\otimes
{\frak g}(x),
$$
i. e.
$$(\id \otimes \Delta){\frak c}+{\frak g}(x)\otimes 1
\otimes 1+1\otimes {\frak c}+1\otimes {\frak g}(x)\otimes 1+
1\otimes 1\otimes {\frak g}(x)=
$$
$$(\Delta \otimes \id ){\frak c}+{\frak c}\otimes 1+
{\frak g}(x)\otimes 1\otimes 1+1\otimes {\frak g}(x)
\otimes 1+1\otimes 1\otimes {\frak g}(x).
$$
This completes the proof.\; $\square$
\begin{remark}
Note that this proof generalizes the standard proof
of the analogous result for formal groups over rings
(see \cite{Honda}).
\end{remark}
\begin{remark}
Note that in the proof we assign for any formal group
${\frak F}(x\otimes 1,1\otimes x)$ over $H$ some 2-cocycle ${\frak c}$
in the cobar complex of the coalgebra $H_\mathbb{Q}.$
\end{remark}
\begin{remark}
Note that $(\varepsilon{\frak g})(x)\in R_\mathbb{Q}[[x]]$
is the logarithm of the formal group $((\varepsilon \otimes
\varepsilon){\frak F})(x\otimes 1,1\otimes x)=F(x\otimes 1,
1\otimes x)\in R[[x\otimes 1,1\otimes x]]$
over ring $R.$
\end{remark}
\begin{remark}
Since ${\frak g}(x)=b_0x+b_1x^2+\ldots\quad \hbox{and}\; \varepsilon(b_0)=1,$
there exists the series
$(\Delta{\frak g})^{-1}(x)=(\Delta({\frak g}^{-1}))(x))\in
H_\mathbb{Q}{\mathop{\widehat{\otimes}}
\limits_{R_\mathbb{Q}}}H_\mathbb{Q}[[x]].$
Using (\ref{4}), we get
$$
{\frak F}(x\otimes 1,1\otimes x)=(\Delta{\frak g})^{-1}
({\frak c}+{\frak g}(x)\otimes 1+1\otimes {\frak g}(x)).
$$
\end{remark}


\begin{thebibliography}{99}
\bibitem{pap1}
{\sc A. V. Ershov}
Floating bundles and their applications.---
arXiv:math.AT/0102054
\bibitem{pap2}
{\sc A. V. Ershov}
Supplement to the paper ''Floating bundles and their
applications''.---
arXiv:math.AT/0102180
\bibitem{Honda}
{\sc Honda T.}
Formal groups and zeta-functions.---
Osaka Journal of Math., 5, 2 (1968), 199--213.
\end{thebibliography}
\end{document}